\documentclass[a4paper,12pt]{article}
\usepackage{mathrsfs}
\usepackage{}
\usepackage[top=2.0cm,bottom=2.0cm,left=2.0cm,right=2.0cm]{geometry}
\usepackage{amssymb}
\usepackage{graphicx}
\usepackage{amsmath,amsthm,amssymb,lineno}
\usepackage{latexsym}
\usepackage{epstopdf}
\usepackage{setspace}
\usepackage{amsmath,color}
\usepackage{graphicx,booktabs,multirow}
\usepackage{latexsym, tabularx,shapepar}
\usepackage[all,2cell,dvips]{xy} \UseAllTwocells \SilentMatrices
\usepackage{appendix}
\usepackage{longtable}
\usepackage{cite}
\usepackage{CJK}
\usepackage{indentfirst}
\usepackage{array}
\usepackage{amsmath}
\usepackage[colorlinks,
           linkcolor=blue,
           anchorcolor=blue,
           citecolor=blue
           ]{hyperref}
\usepackage{dsfont}
\allowdisplaybreaks[4]
\graphicspath{{figures/}}
\usepackage{caption}
\makeatletter
\newcommand{\rmnum}[1]{\romannumeral #1}
\newcommand{\Rmnum}[1]{\expandafter\@slowromancap\romannumeral #1@}
\makeatother

\newtheorem{theorem}{Theorem}[section]
\newtheorem{lemma}[theorem]{Lemma}
\newtheorem{remark}[theorem]{Remark}

\newtheorem{proposition}[theorem]{Proposition}

\newtheorem{claim}{Claim}

\begin{document}

\title{ Maxima of the index: forbidden unbalanced cycles }

\author{Zhuang Xiong\\
	{\small College of Mathematics and Statistics, Hunan Normal University, zhuangxiong@hunnu.edu.cn}\\
	\and
	Yaoping Hou\\
	{\small College of Mathematics and Statistics, Hunan Normal University, yphou@hunnu.edu.cn}}

\maketitle
\begin{abstract}
	This paper aims to address the problem: what is the maximum index among all $\mathcal{C}^-_r$-free unbalanced signed graphs, where $\mathcal{C}^-_r$ is the set of unbalanced cycle of length $r$.
	Let $\Gamma_1 = C_3^- \bullet K_{n-2}$ be a signed graph obtained by identifying a vertex of $K_{n-2}$ with a vertex of $C_3^-$ whose two incident edges in $C_3^-$ are all positive, where $C_3^-$ is an unbalanced triangle with one negative edge. 
	It is shown that if $\Gamma$ is an unbalanced signed graph of order $n$,  $r$ is an integer in  $\{4, \cdots, \lfloor \frac{n}{3}\rfloor + 1 \}$, and 
	$$\lambda_{1}(\Gamma) \geq \lambda_{1}(\Gamma_1), $$ then $\Gamma$ contains an unbalanced cycle of length $r$, unless $\Gamma \sim \Gamma_1$. \\ 
	\indent It is shown that the result are significant in spectral extremal graph problems. 
	Because they can be regarded as a extension of the spectral Tur{\'a}n problem for cycles [Linear Algebra Appl. 428 (2008) 1492--1498] in the context of signed graphs. 
	Furthermore, our result partly resolved a recent open problem raised by Lin and Wang [arXiv preprint arXiv:2309.04101 (2023)].
\end{abstract}
\baselineskip=0.25in

\section{ Introduction}
  Given a signed graph $\Gamma$, the index of $\Gamma$ is defined as the largest eigenvalue of its adjacency matrix $A(\Gamma)$, denoted by $\lambda_1(\Gamma)$.
  In this paper we determine the maximum index of $n$-vertices unbalanced signed graphs with no unbalanced cycle of length $r$, where $r \in \{4, \cdots, \lfloor \frac{n}{3} \rfloor + 1\}$.
  Before presenting our theorem, we will begin with an introductory discussion.\\
  \indent Let $F$ be a set of graphs.  
  A graph $G$ is $F$-free, if it contains no signed graph in $F$ as a subgraph. 
  The spectral Tur{\'a}n problem asks what is the maximum spectral radius of an $F$-free graph of order $n$, which was initially raised by Nikiforovi \cite{nikiforovi2010spectral}.  
  The spectral Tur{\'a}n problem has been studied for many various graphs in the past decades, see $K_r$ \cite{bollobas2007cliques,wilf1986spectral}, $K_{s,t}$ \cite{babai2009spectral,nikiforovi2010contribution}, $C_4$ \cite{nikiforovi2007bounds,zhai2012proof}, $C_5$ \cite{guo2021spectral}, $C_6$ \cite{zhai2020spectral}, consecutive cycles \cite{nikiforovi2008spectral}, and three surveys \cite{chen2018some,nikiforovi2011some,li2022survey}.
  Here we only mention the result of consecutive cycles, which may be stated as follows.
  \begin{theorem}\cite[Theorem 1]{nikiforovi2008spectral}
  	Let $G$ be a graph of sufficiently large order $n$ with $\lambda_{1}(G) > \sqrt{\lfloor n^2/4 \rfloor}$. Then $G$ contains a cycle of length $t$ for every $t \leq n/320$.
  \end{theorem}
  \noindent Recently, Ning and Peng \cite{ning2020extensions} improved above theorem by showing $G$ contains a cycle of length $t$ for every $t \leq n/160$.
In this paper we study a similar problem in the context of signed graphs.\\
  \indent Very recently, Kannan and Pragada first extended the Tur{\'a}n's inequality for signed graphs \cite{kannan2023signed}.
  Afterwards, the spectral Tur{\'a}n problem for signed graphs has been paid much attention. 
  Denote by $\mathcal{K}_k^-$ and $\mathcal{C}_k^-$ the sets of unbalanced signed complete graphs and negative cycles of order $k$, respectively. 
  In 2023, Wang, Hou, and Li \cite{wang2022extremed} investigated the signed graphs with maximum spectral radius over all $\mathcal{C}_3^-$-free unbalanced signed graphs. 
  After that, Chen and Yuan \cite{chen2023turan} studied the spectral Tur{\'a}n problem for $\mathcal{K}_4^-$-free signed graphs.
  Wang \cite{wang2023spectral} continued this study for $\mathcal{K}_5^-$-free signed graph.
  In 2023, Wang and Lin \cite{wang2023largest} considered the existence of a negative $C_{2k}$ from the largest eigenvalue condition, and they solved this problem when $k = 2$ (see Theorem \ref{tho:c4}).
 
\begin{figure}[htbp]\centering\hspace{0cm}
\scalebox{0.5}{\includegraphics[width=8cm]{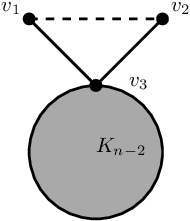}}\\
\caption{ The signed graph $\Gamma_1 = C_3^- \bullet K_{n-2}$. } 
\label{fig:Gamma_1}
\end{figure}

  \begin{theorem}\cite[Theorem 3]{wang2023turan}\label{tho:2k+1cycle}
  	Let $\Gamma$ be an unbalanced signed graph of order $n$ and $3 \leq k \leq \frac{1}{15}n$ be an integer. 
  	If $\lambda_1(\Gamma) \geq \lambda_1(\Gamma_1)$, then $\Gamma$ contains a negative $C_{2k+1}$ unless $\Gamma$ is switching equivalent to $\Gamma_1$ (see Fig. \ref{fig:Gamma_1}).
  \end{theorem}
  
    \begin{theorem}\label{tho:c4}\cite[Theorem 1]{wang2023largest}
  	Let $\Gamma = (G, \sigma)$ be an unbalanced signed graph of order $n \geq 5$. 
  	If $\lambda_{1}(\Gamma) \geq \lambda_1(\Gamma_1)$, then $\Gamma$ contains a negative $C_4$ unless $\Gamma$ is switching equivalent to $\Gamma_1$ (see Fig. \ref{fig:Gamma_1}).
  \end{theorem}
  
  In this work we put our perspective on the existence of a negative cycle $C_r$ from the largest eigenvalue condition. 
  Let $C_3^-$ be a negative triangle with one negative edge. Denote by $\Gamma_1 = C_3^- \bullet K_{n-2}$ the signed graph obtained by identifying a vertex of $K_{n-2}$ with a vertex of $C_3^-$ whose two adjacent edges in $C_3^-$ are all positive (see Fig. \ref{fig:Gamma_1}). 
  Here positive edges (resp., negative edges) are represented by solid lines (resp., dashed lines) and the all positive complete graph is represented by the solid circle.
  
  \begin{theorem}\label{tho:c_r}
  	Let $\Gamma$ be an unbalanced signed graph of order $n$, and $r$ be an integer in $\{4, \cdots, \lfloor \frac{n}{3} \rfloor + 1\}$.
  	If $\lambda_1(\Gamma) \geq \lambda_{1}(\Gamma_1)$, then $\Gamma$ contains a negative cycle of length $r$, unless $\Gamma$ is switching equivalent to $\Gamma_1$.
  \end{theorem}
  
   Note that Wang, Hou, and Huang \cite{wang2023turan} gave the largest eigenvalue condition for the existence of a negative $C_{2k+1}$ when $3 \leq k \leq \frac{1}{15}n$. 
   So Theorem \ref{tho:c_r} improve their result.

  \indent

\section{Preliminaries}\label{sec:pre}
  We need introduce some notations. 
  For signed graph notation and concepts undefined here, we refer the reader to \cite{zaslavskyi1982signed}.
  For spectral theory of signed graphs see the survey \cite{zaslavskyi2010matrices} of Zaslavskyi and the references therein.
  Let $\Gamma$ be an $n$-vertices signed graph, and $X$ be a set of vertices of $\Gamma$.
  We write:\\
  \indent - $V(\Gamma)$ for the set of vertices of $\Gamma$, and $e(\Gamma)$ for the number of its edges; \\
  \indent - $\Gamma[X]$ for the signed graph induced by $X$, and $e(X)$ for $e(\Gamma[X])$;\\
  \indent -$N_{\Gamma}(v)$ for the set of neighbors of a vertex $v$, and $d_{\Gamma}(v)$ for $|N_{\Gamma}(v)|$;\\
  \indent -$u \stackrel{+}{\sim} v$, $u \stackrel{-}{\sim} v$,
  and $u \not \sim v$ for the positive edge, negative edge, and non-edge $uv$, respect;\\
  \indent -$\Gamma^\prime \sim \Gamma$ for that $\Gamma^\prime$ is switching equivalent to $\Gamma$.\\
 
  The following lemma is a key tool for checking the switching equivalence between two signed graphs with the same underlying graph.
  \begin{lemma}\cite[Proposition 3.2]{zaslavskyi1982signed}\label{lem:swi}
  	Two signed graphs on the same underlying graph are switching  equivalent if and only if they have the same list of balanced cycles. 
  \end{lemma}

  In the following two lemmas, we introduce two upper bounds of indexes of a graph and a signed graph, respectively. 
	  The first one was announced by Schwenk in \cite{schwenk1975spectral} and a proof by Hong \cite{yuan1988bound} showed some thirteen years later. 
  The second one posed by Stani{\'c} in \cite{stanic2019bounding} which can be regarded as a signed graph version of the former.
  
  \begin{lemma}\cite[Theorem 1]{yuan1988bound}
  	For any connected graph $G$ with order $n$ and size $e(G)$, $\lambda_{1}(G) \leq \sqrt{2e(G) - n + 1}$.
  \end{lemma}

  Denote by $\epsilon(\Gamma)$ the frustration index of a signed graph $\Gamma$, which refers to the minimum number of edges removed to make the signed graph balanced. 
  The lemma below give a upper bound of the index of a signed graph in terms of its size, order, and frustration index. 
  \begin{lemma}\cite[Theorem 3.2]{stanic2019bounding}\label{lem:fru}
  	Let $\Gamma$ be a connected signed graph with $n$ vertices, $e(\Gamma)$
  	edges, and the frustration index $\epsilon(\Gamma)$. 
  	Then
  \begin{equation*}
  	\lambda_1(\Gamma) \leq \sqrt{2(e(\Gamma)-\epsilon(\Gamma)) - n +1}.
  \end{equation*}
  \end{lemma}

\section{ Proof of Theorem \ref{tho:c_r}}\label{sec:pro}
  Several lemmas which are the ingredients of the proof of Theorem \ref{tho:c_r} will show below. 
  First we give a sharp lower bound for $\lambda_1(\Gamma_1)$.

\begin{lemma}\cite[Lemma 2]{wang2023largest}\label{lem:g1}
	The largest eigenvalue of $\Gamma_1$ is $\lambda_1(\Gamma_1) > n - 3$.
\end{lemma}	

The lemma below show an interesting relationship between the index $\lambda_{1}(\Gamma)$ and its eigenvector in a signed graph $\Gamma$. 
\begin{lemma}\label{lem:zer}\cite[Lemma 3.3]{xiong2023extremal}
	Let $\Gamma = (G, \sigma)$ be a signed graph with a unit eigenvector $x = (x_1, x_2, \cdots, x_n)^{\intercal}$ corresponding to $\lambda_{1}(\Gamma)$. 
	If $\lambda_1(\Gamma) > n - k$, then $x$ has at most $k - 2$ zero component.
\end{lemma}

The following lemma provide a useful technique for studying some perturbations in a signed graph.
\begin{lemma}\cite[Lemma 1]{stanic2018perturbations}\label{lem:nonneg}
	Let $\Gamma$ be a signed graph. 
	Then there exists a signed graph $\Gamma^\prime$ switching equivalent to $\Gamma$ such that $\lambda_1(\Gamma^\prime)$ has a non-negative eigenvector.
\end{lemma}

Let $x = (x_1, x_2, \cdots, x_n)^{\intercal}$ be an eigenvector corresponding to the largest eigenvalue $\lambda_{1}(\Gamma)$ of a signed
graph $\Gamma$. 
The entry $x_i$ is usually corresponding to the vertex $v_i$ of $\Gamma$. 
So the eigenvalue equation for $v_i$ reads as follows
\begin{align*}
	\lambda_{1}(\Gamma)x_{i}  = & \sum_{v_j \in N_{\Gamma}(v_i)} \sigma(v_iv_j)x_j. 
\end{align*}
The following lemma can be proved based on \cite[Theorem 3]{stanic2018perturbations} or \cite[Proposition 2.1]{akbari2019largest}, and a detailed proof can be seen from
\cite[Proposition 3.2]{xiong2023extremal}.
\begin{lemma}\label{lem:per}
	Let $\Gamma = (G, \sigma)$ be a signed graph with a non-negative unit eigenvector $x = (x_1, x_2, \cdots, x_n)^{\intercal}$ corresponding to the largest eigenvalue $\lambda_{1}(\Gamma)$. 
	If we perform one of the following perturbations in $\Gamma$:\\
	\indent \textbf{(\rmnum{1})} Adding some positive edges,\\
	\indent \textbf{(\rmnum{2})} Removing some negative edges,\\
	\indent \textbf{(\rmnum{3})} Reversing the signs of some negative edges,\\
	resulting in a new signed graph $\Gamma ^ \prime $, then $\lambda_{1}(\Gamma^\prime) \geq \lambda_{1}(\Gamma)$. 
	The equality holds if and only if the entries of $x$  corresponding to the endpoints of these edges are all zeros.
\end{lemma}

\begin{remark}
	In view of Lemmas \ref{lem:g1}, \ref{lem:zer}, \ref{lem:nonneg} and \ref{lem:per}, if $\Gamma$ is a $\mathcal{C}_r^-$-free unbalanced signed graph with maximum index, then $\Gamma$ must be connected.
\end{remark}

To simplify the proof of Theorem \ref{tho:c_r}, we shall prove an auxiliary statements which is a crucial case of it.
\begin{proposition}\label{pro:c5}
	Let $\Gamma = (G, \sigma)$ be an unbalanced signed graph of order $n \geq 12$. 
	If $\lambda_{1}(\Gamma) \geq \lambda_1(\Gamma_1)$, then $\Gamma$ contains a negative pentagon unless $\Gamma \sim \Gamma_1$.
\end{proposition}
\begin{proof}
	Suppose that $\Gamma$ has maximum index over all $\mathcal{C}_5^-$-free unbalanced signed graph of order $n \geq 12$.
	Note that from Lemma \ref{lem:g1} we have $\lambda_{1}(\Gamma) > n - 3$.
    In what follows, we will derive a series of claims to show that $\Gamma$ is switching equivalent to $\Gamma_1$.
	
	From Lemma \ref{lem:nonneg} we can find a signed graph $\widetilde{\Gamma} \sim \Gamma$ such that the eigenvector corresponding to $\lambda_1(\widetilde{\Gamma})$ is non-negative.
	So by Lemma \ref{lem:swi}, $\widetilde{\Gamma}$ is also a $\mathcal{C}_5^-$-free unbalanced signed graph and $\lambda_{1}(\widetilde{\Gamma}) = \lambda_{1}(\Gamma) > n - 3$.
	Let $x = (x_1, x_2, \cdots, x_n)^\intercal$ be a non-negative unit eigenvector of $\widetilde{\Gamma}$ corresponding to $\lambda_{1}(\widetilde{\Gamma})$, and note that $x$ has at most one zero component by Lemmas \ref{lem:g1} and \ref{lem:zer}.

   Recall that a chord of a cycle $C$ is an edge joining two non-consecutive vertices of $C$. 
   If there exists a chord of $C$, we say that $C$ is a chorded cycle, otherwise $C$ is a chordless cycle.
   Since $\widetilde{\Gamma}$ is unbalanced, there exist at least  one negative edge and at least one negative cycle. 
   Take a negative cycle $\mathcal{C} = v_1v_2 \cdots v_lv_1$ of the shortest length from $\widetilde{\Gamma}$.
   It is clear that $\mathcal{C}$ is a chordless cycle, otherwise there exists a negative cycle shorter than $\mathcal{C}$.
   It follows that $l \leq 4$, otherwise $\widetilde{\Gamma}$ is $\mathcal{C}_4^-$-free and so $\lambda_{1}(\widetilde{\Gamma}) < \lambda_{1}(\Gamma_1)$ by Theorem \ref{tho:c4}, a contraction.

	\begin{claim}\label{cla:loc}
		All the negative edges of $\widetilde{\Gamma}$ are lied in $\mathcal{C}$.
		\begin{proof}
			If not, assume that there is a negative edge not in $\mathcal{C}$.
			Then we can delete it from $\widetilde{\Gamma}$ to obtain a $\mathcal{C}_5^-$-free unbalanced signed graph $\widetilde{\Gamma}^\prime$.
			By Lemma \ref{lem:per}, we have $\lambda_1({\widetilde{\Gamma}^\prime}) > \lambda_1({\widetilde{\Gamma}})$, a contraction.		
		\end{proof}
	\end{claim}
	
	\begin{claim}\label{cla:one}
		$\widetilde{\Gamma}$ contains only one negative edge. 
	\end{claim}
		\begin{proof}
			If not, from Claim \ref{cla:loc} and $l \leq 4$, assume that there exist $3$ negative edges which are contained in $\mathcal{C}$.\\
			\indent For $l = 3$, we claim that there is no positive path $P$ together with two negative edges of $\mathcal{C}$ forming a positive $C_5$. 
			Indeed, without loss of generality, if $P$, $v_1v_2$, and $v_2v_3$ form a positive $C_5$, then by Lemma \ref{lem:per} we may delete $v_1v_2$ and $v_2v_3$ we obtain a $\mathcal{C}_5^-$-free unbalanced signed graph with lager index than $\widetilde{\Gamma}$, a contradiction. 
			Hence, by reversing the signs of any two negative edges of $\mathcal{C}$, we also obtain a $\mathcal{C}_5^-$-free unbalanced signed graph with lager index than $\widetilde{\Gamma}$, a contradiction.\\
			\indent For $l = 4$, first note that there is no positive path $P$ together with two consecutive negative edges forming a positive $C_5$, otherwise $P$ and remaining two edges form a negative $C_5$, a contradiction.
			Furthermore, there exist no positive path $P$ together with two nonconsecutive negative edges forming a positive $C_5$, otherwise there exist a negative triangle, a contradiction to $l = 4$.
			Hence, by Lemma \ref{lem:per} we may reverse any two negative edges of $\mathcal{C}$, we also obtain a $\mathcal{C}^-_5$-free unbalanced signed graph with larger index than $\widetilde{\Gamma}$, a contradiction.
		\end{proof}

	From Claim \ref{cla:one} we have $\epsilon(\widetilde{\Gamma}) = 1$.
	Without loss of generality, we always  assume that $v_1v_2$ is that negative edge and $e(\widetilde{\Gamma}) = \frac{n(n-1)}{2} - h$.
	In what follows, we first obtain a tight upper bound for $h$ which is a crucial quantity for characterizing the structure of $\widetilde{\Gamma}$.
	
	\begin{claim}\label{cla:number}
		$h \leq 2n - 6$.
		\begin{proof}
			If not, $h \geq 2n -5$. By Lemma \ref{lem:fru}, we have
			\begin{align*}
				\lambda_{1}(\widetilde{\Gamma}) &\leq \sqrt{2(e(\widetilde{\Gamma}) - \epsilon(\widetilde{\Gamma})) - n + 1} \\
				&\leq \sqrt{2(\frac{n(n-1)}{2} - 2n + 5 - 1) - n + 1}\\
				&=\sqrt{n^2 - 6n + 9}\\
				&=n-3,
			\end{align*}
			which contradict to Lemma \ref{lem:g1}.
		\end{proof}
	\end{claim}

	\begin{claim}\label{cla:length}
		$l = 3$.
	\end{claim}
		\begin{proof}
		If not, then $l = 4$.
		Bear in mind that $\mathcal{C}$ is a chordless cycle.\\
		\indent We assert that $v_i$ is adjacent to at most two vertices of $V(\mathcal{C}) = \{v_1,v_2,v_3,v_4 \}$, for any vertex $v_i \in V(\widetilde{\Gamma}) \setminus V(\mathcal{C})$.
		Indeed, if $\{v_1,v_4\} \subseteq N_{\widetilde{\Gamma}}(v_i)$ , then $v_iv_4v_3v_2v_1v_i$ is a negative pentagon, a contradiction. 
		Similarly, $\{v_3, v_4\} \not \subseteq N_{\widetilde{\Gamma}}(v_i)$ and $\{v_2, v_3\} \not \subseteq N_{\widetilde{\Gamma}}(v_i)$.
		Furthermore, $\{v_1, v_2 \} \not \subseteq N_{\widetilde{\Gamma}}(v_i)$, otherwise there is a negative triangle $v_iv_1v_2v_i$.
		So far the number of edges we have lost is at least $ 2 + 2(n-4) = 2n -6$. 
		If there exists a vertex $v_i$ only adjacent to one vertex of $V(\mathcal{C})$, then $h \geq 2n - 5$, which contradicts to Claim \ref{cla:number}.
		Thus, each vertex $v_i \in V(\widetilde{\Gamma}) \setminus V(\mathcal{C})$ is exactly adjacent to two vertices of $V(\mathcal{C})$, that is, $\{v_1, v_3 \}$ or $\{v_2, v_4 \}$ $\subseteq N_{\widetilde{\Gamma}}(v_i)$.
		We assert that both $N_{\widetilde{\Gamma}}(v_1) \cap N_{\widetilde{\Gamma}}(v_3)
		\setminus \{v_2,v_4\}$  and $N_{\widetilde{\Gamma}}(v_2) \cap N_{\widetilde{\Gamma}}(v_4) \setminus \{v_1,v_3\}$ are independent sets. 
		Here we prove the former, and the latter is completely similar to it.
		In fact,  if there exist two vertices $v_i, v_j$ such that $\{v_i,v_j\} \subseteq N_{\widetilde{\Gamma}}(v_1) \cap N_{\widetilde{\Gamma}}(v_3) \setminus \{v_2,v_4\}$, then $v_i \nsim v_j$, otherwise $v_iv_jv_3v_2v_1v_i$ is a negative pentagon, a contradiction.
		Since $|V(\widetilde{\Gamma}) \setminus V(\mathcal{C}) | \geq 8$, then we lose at least one edge again and so $h \geq 2n -5$, a contradiction.
		\end{proof}
	
	From Claim \ref{cla:length} we know that $\mathcal{C}$ is a negative triangle, say $v_1v_2v_3v_1$, with the negative edge $v_1v_2$.
	Recall that $x_i$ is the entry of $x$ corresponding to the vertex $v_i$.
	Let $x_m = \max \{x_1, \cdots, x_n\}$.
	
	\begin{claim}\label{cla:n-1}
		$d_{\widetilde{\Gamma}}(v_m) = n - 1$.
	\end{claim}
		\begin{proof}
			If not, $d_{\widetilde{\Gamma}}(v_m) \leq n - 2$.\\
			\indent First we consider $d_{\widetilde{\Gamma}}(v_m) \leq n - 3$.
			By the eigenvalue equation for $v_m$, we have 
			\begin{align*}
				\lambda_{1}(\widetilde{\Gamma})x_m = \sum_{v_k  \in N_{\widetilde{\Gamma}}(v_m)} \sigma(v_kv_m)x_k \leq (n-3)x_m,
			\end{align*}
			and so $\lambda_{1}(\widetilde{\Gamma}) \leq n - 3$, a contradiction.\\
			\indent Then we have $d_{\widetilde{\Gamma}}(v_m) = n - 2$.
			Next we will divide our proof into three cases in terms of the value of $m$.\\
			\indent \textbf{Case 1.} $m = 1$.\\
			\indent Assume without loss of generality that $N_{\widetilde{\Gamma}}(v_m) = \{v_2, \cdots v_{n-1} \}$. 
			Then
				\begin{align*}
					\lambda_{1}(\widetilde{\Gamma})x_1 = - x_2 +  \sum_{3 \leq i \leq n - 1} x_i \leq (n-3)x_1,
				\end{align*}
				and so $\lambda_{1}(\widetilde{\Gamma}) \leq n - 3$, a contradiction.
				The case of $m = 2$ is the same with the above. \\
			\indent \textbf{Case 2.} $m = 3$.\\
			\indent Without loss of generality, suppose that $v_3 \nsim v_n$.\\
			\indent \textbf{Subcase 1.}
				$v_n \sim v_1$, see Fig. \ref{fig:r3ab} $(a)$. \\ 
				We claim that $v_4, v_5, \cdots, v_{n-1} \notin N_{\widetilde{\Gamma}}(v_n)$, otherwise there is a vertex, say $v_4$, adjacent to $v_n$, then $v_nv_1v_2v_3v_4v_n$ is a negative pentagon (shown in Fig. \ref{fig:r3ab} $(a_1)$), a contradiction.
				If there exists an edge between $\{v_1\}$ and $\{v_4, \cdots, v_{n-1} \}$, say $v_1v_4$ , 
				then $ v_5, \cdots, v_{n-1} \notin N_{\widetilde{\Gamma}}(v_2)$, otherwise there is a vertex in $\{v_5, \cdots, v_{n-1} \}$ adjacent to $v_2$, say $v_5$, such that $v_2v_1v_4v_3v_5v_2$ is a negative pentagon (shown in Fig. \ref{fig:r3ab} $(a_2)$), a contradiction. 
				Meanwhile, $ v_5, \cdots, v_{n-1} \notin N_{\widetilde{\Gamma}}(v_4)$, otherwise there exists a negative pentagon $v_2v_1v_4v_iv_3v_2$ for $i \in \{5, \cdots, n-1\}$ (shown in Fig. \ref{fig:r3ab} $(a_3)$), a contradiction. 
				So far, the number of edges we have lost is at least $1 + n - 4 + 2(n - 5) = 3n - 13 > 2n - 6$, which contradict to Claim \ref{cla:number}.
				Therefore, $v_4, \cdots, v_{n-1} \notin N_{\widetilde{\Gamma}}(v_1)$. 
				By the similar analysis we have $v_4, \cdots, v_{n-1} \notin N_{\widetilde{\Gamma}}(v_2)$.
				So the number of edges we have lost is at least $1 + n - 4 + 2(n - 4)  = 3n - 11 > 2n - 6 $, a contradiction.
				The case of $v_n \sim v_2$ is completely similar to the above. \\
			\indent \textbf{Subcase 2.} $v_n \nsim v_1, v_2$, see Fig. \ref{fig:r3ab} $(b)$. \\
				If there exists an edge between $\{v_1\}$ and $\{v_4, \cdots, v_{n-1} \}$, say $v_1v_4$,  by the same discussion to the above, $ v_5, \cdots, v_{n-1} \notin N_{\widetilde{\Gamma}}(v_2) \cup N_{\widetilde{\Gamma}}(v_4)$.
				So we have lost $3 + 2(n - 5) = 2n - 7$ edges and at most one edge we can lose now. 
			    Then there must exist an edge $v_iv_j$ such that $v_2v_1v_iv_jv_3v_2$ is a negative pentagon (shown in Fig. \ref{fig:r3ab} $(b_1)$), a contradiction. 
			    Thus, $v_4, \cdots, v_{n-1} \notin N_{\widetilde{\Gamma}}(v_1)$ and by symmetry we also have $v_4, \cdots, v_{n-1} \notin N_{\widetilde{\Gamma}}(v_2)$.
				Then the number of edges we have lost is at least $3 + 2(n - 4) = 2n - 5$, a contradiction.\\
			\indent \textbf{Case 3.} $m \geq 4$.\\ 
			\indent Without loss of generality assume that $m = 4$.	
			Then there exists a vertex being not adjacent to $v_4$.		
				If $v_4 \nsim v_2$, shown in Fig. \ref{fig:r3ab} $(c)$, then $ v_5, \cdots, v_{n} \notin N_{\widetilde{\Gamma}}(v_1) \cup N_{\widetilde{\Gamma}}(v_2) \cup N_{\widetilde{\Gamma}}(v_3)$.
				So we lose $1 + 3(n - 4) = 3n - 11 > 2n -6$ edges, a contraction.
				The discussion of $v_4 \nsim v_1$ is almost the same with the above.\\
				\indent If $v_4 \nsim v_3$, shown in Fig. \ref{fig:r3ab} $(c_1)$, then $ v_5, \cdots, v_{n} \notin  N_{\widetilde{\Gamma}}(v_3)$.
				Further, if there is a vertex $v_i \in \{v_5, \cdots, v_{n}\}$ adjacent to $v_1$, say $v_5$,  then $ v_6, \cdots, v_{n} \notin  N_{\widetilde{\Gamma}}(v_2) \cup N_{\widetilde{\Gamma}}(v_5)$ and we lose $1 + n - 4 + 2(n - 5) = 3n - 13 > 2n - 6$ edges, a contradiction.
				So $ v_5, \cdots, v_{n} \notin   N_{\widetilde{\Gamma}}(v_1)$ and by symmetry we have $ v_5, \cdots, v_{n} \notin   N_{\widetilde{\Gamma}}(v_2)$. 
				We have lost $3(n - 4) + 1 = 3n -11 > 2n - 6$ edges, a contradiction.\\
				\indent If $v_4$ is not adjacent to a vertex of $\{v_5, \cdots, v_n \}$, say $v_n$, shown in Fig. \ref{fig:r3ab} $(c_2)$, then $ v_5, \cdots, v_{n-1	} \notin  N_{\widetilde{\Gamma}}(v_1) \cup N_{\widetilde{\Gamma}}(v_2) \cup N_{\widetilde{\Gamma}}(v_3)$, and we have lost $3(n - 5) + 1 = 3n - 14 > 2n -6$ edges, a contradiction.\\
				\indent Now we have considered all the cases, and so Claim \ref{cla:n-1} holds.
					\begin{figure}[htbp]\centering\hspace{0cm}
					\scalebox{1}{\includegraphics[width=16cm]{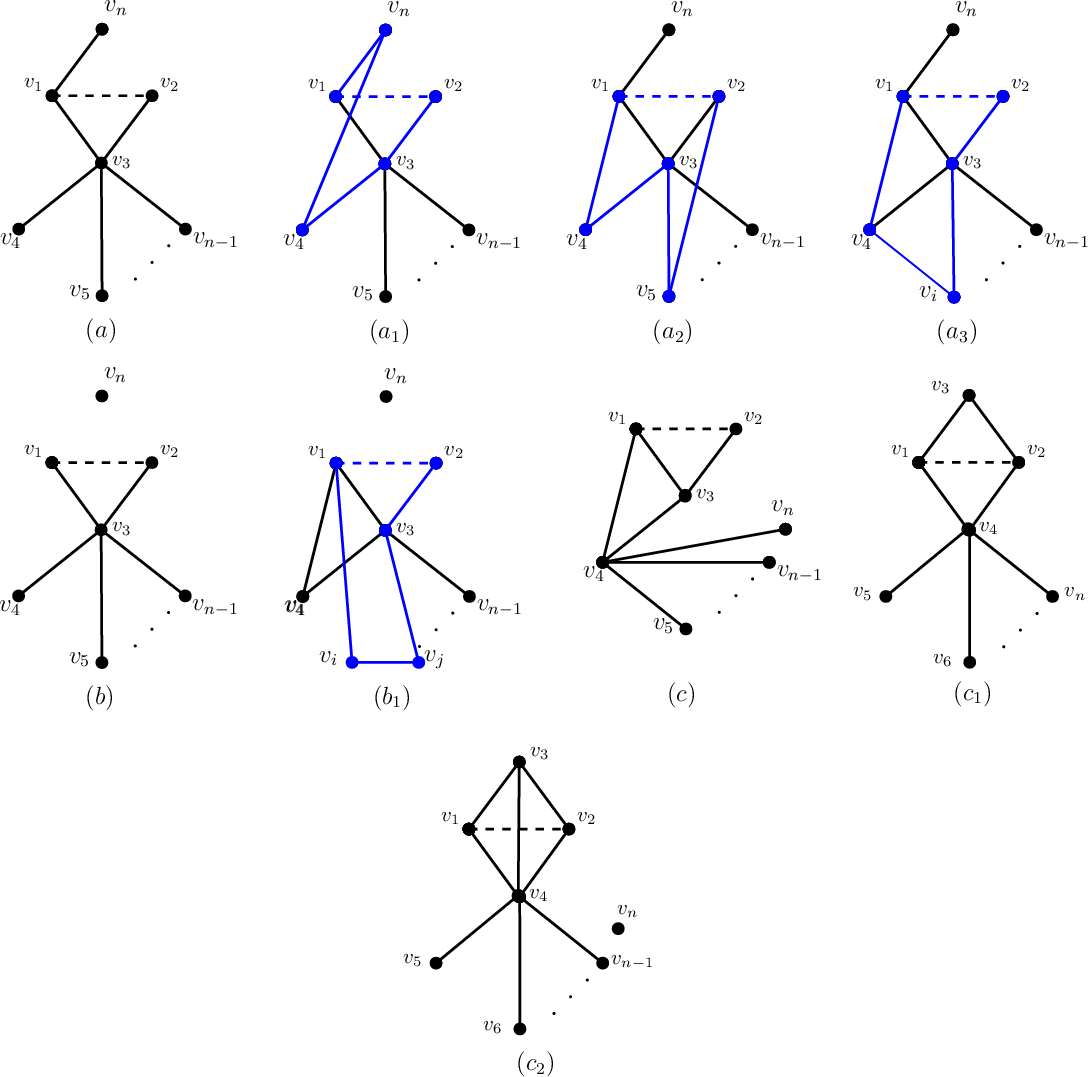}}\\
					\caption{ The signed graphs mentioned in the proof of Claim \ref{cla:n-1}. } 
					\label{fig:r3ab}
					\end{figure}
	\end{proof}
				
		\begin{claim}\label{cla:m3}
		m = 3.
	\end{claim}
	\begin{proof}
		If not, we divide the proof into two cases in terms of the value  of $m$.\\
		\indent \textbf{Case 1.} $m = 1$.\\
		\indent \textbf{Subcase 1.} There is an edge between $\{v_2\}$ and $\{v_4, \cdots, v_n  \}$, say $v_2v_4$. \\
		\indent If $v_3 \sim v_4$, then $ v_5, \cdots, v_{n} \notin  N_{\widetilde{\Gamma}}(v_3) \cup N_{\widetilde{\Gamma}}(v_4).$
		So we have lost $2(n - 4) = 2n - 8$ edges and at most two edges we can lose later.
		Then there must exist two vertices in $ \{v_6, \cdots, v_{n}\}$, say $v_6$ and $v_n$, such that $v_2v_1v_5v_6v_nv_2$ is a negative pentagon, shown in Fig. \ref{fig:m3ab} $(a)$, a contradiction.\\
		\indent So $v_3 \nsim v_4$ and note that each vertex in $ \{  v_5, \cdots, v_{n} \}$ can not be adjacent to both $v_3$ and $v_4$, otherwise say $v_5 \sim v_3$ and $v_5 \sim v_4$, then $v_2v_1v_3v_5v_4v_2$ is a negative pentagon, shown in 
		Fig. \ref{fig:m3ab} $(a_1)$, a contradiction. We lose at least $n - 4$ edges here.
		If $ v_5, \cdots, v_{n} \notin  N_{\widetilde{\Gamma}}(v_3) \cup N_{\widetilde{\Gamma}}(v_4)$, then we can also find a negative pentagon like the above, a contradiction. \\
		\indent If there is a vertex in $\{v_5, \cdots, v_{n} \}$ adjacent to $v_4$, say $v_5$, then $v_6, \cdots, v_n \notin N_{\widetilde{\Gamma}}(v_5)$.
		Further, if there is a vertex in $\{v_6, \cdots, v_{n} \}$ adjacent to $v_3$, say $v_n$, then $v_6, \cdots, v_{n-1} \notin N_{\widetilde{\Gamma}}(v_n)$, shown in Fig.  \ref{fig:m3ab} $(a_2)$.
		So we will have lost  at least $1 + (n - 4) + (n - 5) + (n -6) = 3n - 14 > 2n - 6$ edges, a contradiction.
		Thus $v_5, \cdots, v_{n} \notin  N_{\widetilde{\Gamma}}(v_3)$. 
		If there exists another vertex in $\{v_6, \cdots, v_{n} \}$ adjacent to $v_4$, say $v_6$, then $v_7, \cdots, v_n \notin N_{\widetilde{\Gamma}}(v_6)$, and we lose at least $1 + n - 4 + n - 5 + n -6 = 3n - 14 > 2n -6$ edges, a contradiction. 
		So $\{v_6, \cdots, v_{n} \} \notin N_{\widetilde{\Gamma}}(v_4)$, and we have lost at least $1 + n - 4 + n - 5 + n - 5 = 3n - 13 > 2n -6$ edges, a contradiction.
		Thus, $ v_5, \cdots, v_{n} \notin  N_{\widetilde{\Gamma}}(v_4)$, and by symmetry $ v_5, \cdots, v_{n} \notin  N_{\widetilde{\Gamma}}(v_3)$, we have a contradiction like above. \\
		\indent \textbf{Subcase 2.} There is no edge between $\{v_2\}$ and $\{v_4, \cdots, v_n \}$.\\
		\indent If there is an edge between $\{v_3\}$ and $\{v_4, \cdots, v_n\}$, say $v_3v_n$, shown in Fig. \ref{fig:m3ab} $(b)$, then $ v_4, \cdots, v_{n-1} \notin  N_{\widetilde{\Gamma}}(v_n)$.
		So we have lost $n - 3 + n - 4 = 2n - 7$ edges, and there must exist a negative pentagon $v_2v_1v_4v_iv_3v_2$ for $v_i
		\in \{v_5, \cdots v_{n-1}\}$, shown in Fig. \ref{fig:m3ab} $(b_1)$, a contradiction.
		We have known that there is no edge between $\{v_2, v_3\}$ and $\{v_4, \cdots, v_n\}$.
		Thus, 
			\begin{align*}
				\lambda_{1}(\widetilde{\Gamma})x_2 &= -x_1 + x_3,\\
				\lambda_{1}(\widetilde{\Gamma})x_3 &= x_1 + x_2,
			\end{align*}
			and so $\lambda_{1}(\widetilde{\Gamma})(x_2 + x_3) = x_2 + x_3$, a contradiction.		
		The proof of $m = 2$ is the same with the above.\\
		\indent \textbf{Case 2.} $m \geq 4$. \\
		\indent Without loss of generality, assume that $m = 4$, shown in Fig. \ref{fig:m3ab} (c).  Then $v_5, \cdots, v_n \notin N_{\widetilde{\Gamma}}(v_1) \cup N_{\widetilde{\Gamma}}(v_2) \cup N_{\widetilde{\Gamma}}(v_3)$. We have lost at least $3(n-4) = 3n - 12 > 2n - 6$ edges, a contradiction.
		So the claim holds.
	\end{proof}
	
	\begin{claim}\label{cla:finally}
		$\widetilde{\Gamma} = \Gamma_1$ 
	\end{claim}
	\begin{proof}
		From Claim \ref{cla:m3} we have $d_{\widetilde{\Gamma}}(v_3) = n - 1$, shown in Fig. \ref{fig:finally} $(a)$.
		If there is a vertex in $\{v_4, \cdots, v_n\}$, say $v_4$, adjacent to both $v_1$ and $v_2$, shown in Fig. \ref{fig:finally} $(a_1)$, then $ v_5, \cdots, v_n \notin
		N_{\widetilde{\Gamma}}(v_1) \cup
		N_{\widetilde{\Gamma}}(v_2) \cup N_{\widetilde{\Gamma}}(v_4)$.
		So we lose $3(n-4) = 3n - 12 > 2n -6$ edges, a contradiction.\\
		\indent So if there is an edge between $ \{v_1, v_2 \}$ and $\{v_4, \cdots, v_n\}$, say $v_1v_4$, shown in Fig. \ref{fig:finally} $(a_2)$, then $ v_5, \cdots, v_n \notin
		N_{\widetilde{\Gamma}}(v_2) \cup
		N_{\widetilde{\Gamma}}(v_4)$ and $v_2 \nsim v_4$.
		We have lost $2(n-4) + 1 = 2n - 7$ edges, and we can lose at most one more edge. 
		It is easy to see the existence of a negative pentagon $v_2v_1v_iv_jv_3v_2$ for $v_i,v_j \in \{  v_5, \cdots, v_n\}$, a contradiction.
		Therefore, there is no edge between $ \{v_1, v_2 \}$ and $\{v_4, \cdots, v_n\}$ and we have lost $2n - 6$ edges.
		Then $\{v_3, v_4, \cdots, v_n\}$ is a clique, that is, $\widetilde{\Gamma} = \Gamma_1$.
	\end{proof}
	
		\begin{figure}[htbp]\centering\hspace{0cm}
			\scalebox{1}{\includegraphics[width=16cm]{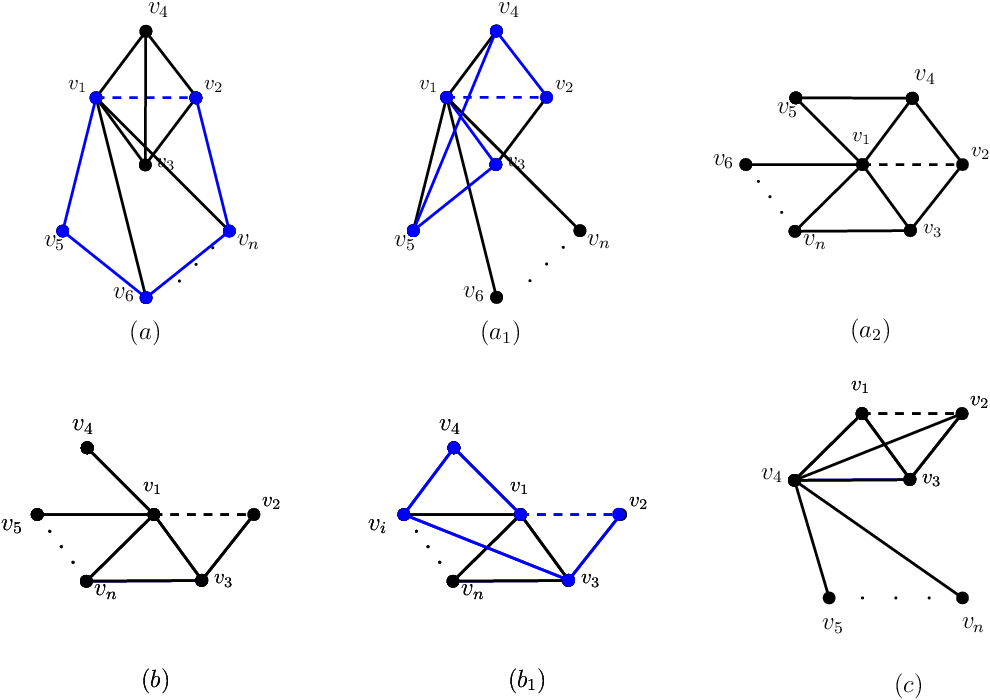}}\\
			\caption{ The signed graphs mentioned in the proof of Claim \ref{cla:m3}. } 
			\label{fig:m3ab}
		\end{figure}
	
		\begin{figure}[htbp]\centering\hspace{0cm}
		\scalebox{1}{\includegraphics[width=16cm]{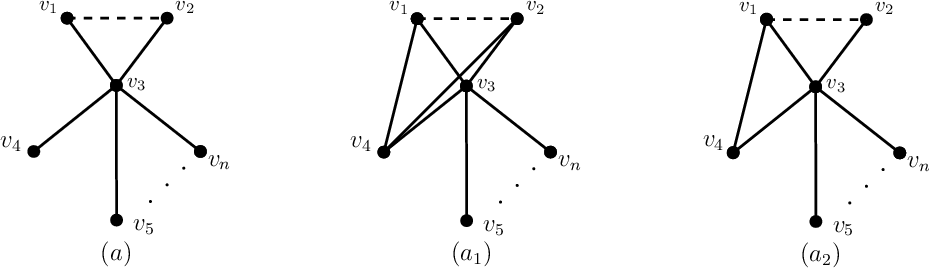}}\\
		\caption{ The signed graphs mentioned in the proof of Claim \ref{cla:finally}. } 
		\label{fig:finally}
	\end{figure}
	In summary, $\Gamma \sim \widetilde{\Gamma} = \Gamma_1$,  completing the proof of Proposition \ref{pro:c5}.
\end{proof}	

\begin{figure}[htbp]\centering\hspace{0cm}
	\scalebox{0.5}{\includegraphics[width=8cm]{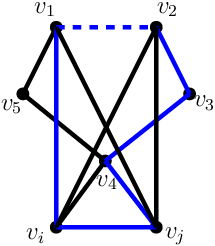}}\\
	\caption{ The signed graph containing a negative hexagon. } 
	\label{fig:hexagon}
\end{figure}

\begin{proof}[Proof of Theorem \ref{tho:c_r}]
	The cases of $r = 4, 5$ have been proven in Theorem \ref{tho:c4} and Proposition \ref{pro:c5}, respectively.
	The following proof relies on induction hypothesis for $r$.
	Assume that the result holds for $\mathcal{C}^-_i$-free unbalanced graphs, $4 \leq i \leq r -1$.
	Next we consider the $\mathcal{C}^-_r$-free unbalanced signed graph $\Gamma$ and let $\Gamma \sim \widetilde{\Gamma} $ which has a non-negative eigenvector.\\
	\indent If $\widetilde{\Gamma}$ does not contain a negative $C_i$ for some $i \in \{4, \cdots, r - 1 \}$, by induction hypothesis we conclude $\widetilde{\Gamma} \sim \Gamma_1$. \\
	\indent Therefore, $\widetilde{\Gamma}$ contains negative cycles of length from $4$ to $r-1$.
	By the same discussion as Claim \ref{cla:one}, there exists a negative $C_{r - 1}$ with one negative edge $v_1v_2$, say $C^-_{r - 1}$. 
	We claim that the number of edges between $V(C^-_{r - 1})$ and $V(\widetilde{\Gamma}) \setminus V(C^-_{r - 1})$ are at most $\lceil \frac{r - 1}{2} \rceil (n - r + 1)$.
	Indeed, each $v_i \in V(\widetilde{\Gamma}) \setminus V(C^-_{r - 1})$ can not be adjacent to two consecutive vertices other than $v_1$ and $v_2$, otherwise there is a negative $C_r$.
	On the one hand, 
	$$ e(\widetilde{\Gamma}) \leq \frac{n(n-1)}{2} - \lfloor \frac{r - 1}{2} \rfloor (n - r + 1). $$
	One the other hand, we have 
	$$ e(\widetilde{\Gamma}) \geq \frac{n (n - 1)}{2} - 2(n - 3).$$
	Note that, for $r \geq 7$, 
	$$ \lfloor \frac{r - 1}{2} \rfloor (n - r + 1) - 2(n - 3) \geq 3(n - r + 1) - 2n + 6  = n - 3r + 9 > 0,$$
	a contradiction.
	Next we consider the case of $r = 6$ and write $v_1v_2v_3v_4v_5v_1$ for the $C^-_5$. We have
	$$ e(\widetilde{\Gamma}) \leq \frac{n(n-1)}{2} -  2 (n - 5), $$
	which means that we can lose at most $4$ more edges.
	If there exist $5$ vertices in $V(\widetilde{\Gamma}) \setminus V(C_5^-)$ being adjacent to at most $2$ vertices of $V(C_5^-)$, then we lost $5$ edges again, a contradiction. 
	Hence, there exist at least $n - 4 - 5 \geq 6$ vertices in $V(\widetilde{\Gamma}) \setminus V(C_5^-)$ being adjacent to exactly $3$ vertices of $V(C_5^-)$.
	It is clear that these vertices must be adjacent to $v_1$, $v_2$, and $v_5$.
	So we claim that these vertices consist of an independent set of $\widetilde{\Gamma}$. Indeed, if $v_i \sim v_j$, then $v_1v_2v_3v_4v_jv_iv_1$ is a negative hexagon (see Fig. \ref{fig:hexagon}), a contradiction. 
	So we lose at least $\frac{(n - 9)(n - 10)}{2} \geq 5$ edges again, a contradiction. This completes the proof of Theorem \ref{tho:c_r}.
\end{proof}

\noindent \textbf{Declaration of competing interest}

The authors declare that there is no conflict of interest.

\vskip 0.6 true cm
\noindent {\textbf{Acknowledgments}}

This research is supported by the National Natural Science Foundation of China (No.11971164).


\baselineskip=0.25in


\begin{thebibliography}{99}

\bibitem{akbari2019largest} S. Akbari, F. Belardo,  F. Heydari, M. Maghasedi, M. Souri, On the largest eigenvalue of signed unicyclic graphs. Linear Algebra Appl. 581 (2019) 145--162.

\bibitem{babai2009spectral} L. Babai, B. Guiduli, Spectral extrema for graphs: the Zarankiewicz problem, Electron. J. Comb. 16 (1) (2009) 123.

\bibitem{bollobas2007cliques} B. Bollob{\'a}, V. Nikiforov, Cliques and the spectral radius, J. Comb. Theory, Ser. B 97 (2007)
859--865.

\bibitem{chen2023turan} F. Chen , X. Y. Yuan, Tur\'{a}n problem for $\mathcal{K}_4^-$-free signed graphs, arXiv preprint arXiv:2306.06655 (2023).

\bibitem{chen2018some} M. Z. Chen, X. D. Zhang, Some new results and problems in spectral extremal graph theory (in Chinese) J. Anhui Univ. Nat. Sci. 42 (2018) 12--25.

\bibitem{guo2021spectral} H. T. Guo, H. Q. Lin,  Y. H. Zhao, A spectral condition for the existence of a pentagon in non-bipartite graphs, Linear Algebra Appl. 627 (2021) 140--149.

\bibitem{yuan1988bound} Y. Hong, A bound on the spectral radius of graphs, Linear Algebra Appl. 108 (1988) 135--139.

\bibitem{kannan2023signed} M. R. Kannan, S. Pragada, Signed spectral Tur{\'a}n type theorems, Linear Algebra Appl. 663 (2023): 62--79.

\bibitem{li2022survey} Y. T. Li, W. J. Liu, L. H. Feng, A survey on spectral conditions for some extremal graph problems, Adv. Math. (China) 51 (2) (2022) 193--258.

\bibitem{nikiforovi2007bounds} V. Nikiforov, Bounds on graph eigenvalues II, Linear Algebra Appl. 427 (2007) 183--189.

\bibitem{nikiforovi2008spectral} V. Nikiforov, A spectral condition for odd cycles in graphs, Linear Algebra Appl. 428 (2008) 1492--1498.

\bibitem{nikiforovi2010contribution} V. Nikiforov, A contribution to the Zarankiewicz problem, Linear Algebra Appl. 432 (2010) 1405--1411.

\bibitem{nikiforovi2010spectral} V. Nikiforov, The spectral radius of graphs without paths and cycles of specified length, Linear Algebra Appl. 432 (2010) 2243--2256.

\bibitem{nikiforovi2011some} V. Nikiforov, Some new results in extremal graph theory, in: Surveys in Combinatorics 2011, in: London Math. Society Lecture Note Ser., vol. 392, 2011, pp. 141--181.

\bibitem{ning2020extensions} B. Ning, X. Peng, Extensions of the Erd\H{o}s-Gallai theorem and Luo's theorem, Combin. Probab. Comput. 29 (2020) 128--136.

\bibitem{schwenk1975spectral} A. J. Schwenk, New derivations of spectral bounds for the chromatic number, Graph Theory Newsletter 5 (1975), No. 1, 77.

\bibitem{stanic2018perturbations} Z. Stani{\'c}, Perturbations in a signed graph and its index, Discuss. Math. Graph T. 38 (3) (2018): 841--852. 

\bibitem{stanic2019bounding} Z. Stani{\'c}, Bounding the largest eigenvalue of signed graphs, Linear Algebra Appl. 573 (2019): 80--89.

\bibitem{wang2022extremed} D. J. Wang, Y. P. Hou, D. Q. Li,  Extremal results for $C_3^-$-free signed graphs, Linear Algebra Appl. 681 (2024): 47--65.

\bibitem{wang2023turan} J. J. Wang, Y. P. Hou, X. Y. Huang, Tur\'{a}n problem for $C_{2k+1}^-$-free signed graph, arXiv preprint arXiv:2310.11061 (2023).

\bibitem{wang2023spectral} Y. A. Wang, Spectral Tur\'{a}n problem for $\mathcal{K}_5^-$-free signed graphs, arXiv preprint arXiv:2309.15434 (2023).

\bibitem{wang2023largest} Y. A. Wang, H. Q. Lin, The largest eigenvalue of $\mathcal {C} _4^{-} $-free signed graphs, arXiv preprint arXiv:2309.04101 (2023).

\bibitem{wilf1986spectral} H. Wilf, Spectral bounds for the clique and independence numbers of graphs, J. Comb. Theory, Ser. B 40 (1986) 113--117.

\bibitem{xiong2023extremal} Z. Xiong, Y. P. Hou, Extremal results for $\mathcal{K}^-_{r + 1}$-free signed graphs.

\bibitem{zaslavskyi1982signed} T. Zaslavsky, Signed graphs, Discrete Appl. Math. 4 (1982) 47--74.

\bibitem{zaslavskyi2010matrices} T. Zaslavsky, Matrices in the theory of signed simple graphs, Advances in Discrete Mathematics and Applications: Mysore, 2008, Ramanujan Math. Soc., Mysore, 2010, pp. 207--229.

\bibitem{zhai2012proof} M. Q. Zhai, B. Wang, Proof of a conjecture on the spectral radius of $C_4$-free graphs, Linear Algebra Appl. 437 (2012) 1641--1647.

\bibitem{zhai2020spectral} M. Q. Zhai, H. Q. Lin, Spectral extrema of graphs: forbidden hexagon, Discrete Math. 343 (10) (2020) 112028.




\end{thebibliography}
\end{document}